\documentclass[11pt]{article}

\usepackage{amsmath}
\usepackage{amssymb}
\usepackage{amsthm}
\usepackage[dvips]{graphicx}
\usepackage{color}
\usepackage{epsfig}
\usepackage[mathscr]{eucal}
\newtheorem{theorem}{Theorem}

\newtheorem{prop}[theorem]{Proposition}
\newtheorem{lemma}[theorem]{Lemma}
\newtheorem{cor}[theorem]{Corollary}
\newtheorem*{claim}{Claim}
\newtheorem{definition}[theorem]{Definition}
\newenvironment{poc}{\noindent \emph{Proof of Claim.}}{\vspace*{0.5\baselineskip}}

\allowdisplaybreaks

\newcommand{\cl}[2]{\ensuremath{\langle {#1} | {#2} \rangle}}

\def\ispreprint{}

\ifx\ispreprint\undefined
\def\baselinestretch{2.5}
\else
\fi

\author{Joshua Cooper\\ \small Department of Mathematics \\ \small University of California, San Diego, La Jolla, CA \\ \small \texttt{jcooper@math.ucsd.edu}}
\date{\today}

\begin{document}
\title{{\huge \bf Quasirandom Permutations}}
\ifx\ispreprint\undefined
\author{{\Large Joshua N.\ Cooper}\\
{\it Department of Mathematics, University of California}\\
{\it at San Diego, La Jolla, California}\\
E-mail: jcooper@math.ucsd.edu\\
Proposed running head: Quasirandom Permutations}
\else
\author{{\Large Joshua N.\ Cooper}\\
{\it Department of Mathematics, University of California}\\
{\it at San Diego, La Jolla, California}\\
E-mail: jcooper@math.ucsd.edu}
\fi

\date{July 1, 2002}

\maketitle

\ifx\ispreprint\undefined \pagebreak \else \fi
\renewcommand{\thesection}{\Roman{section}.}
\maketitle

\begin{abstract}
Chung and Graham [\ref{CG2}] define quasirandom subsets of
$\mathbb{Z}_n$ to be those with any one of a large collection of
equivalent random-like properties.  We weaken their definition and
call a subset of $\mathbb{Z}_n$ $\epsilon$-balanced if its
discrepancy on each interval is bounded by $\epsilon n$.  A
quasirandom permutation, then, is one which maps each interval to
a highly balanced set.  In the spirit of previous studies of
quasirandomness, we exhibit several random-like properties which
are equivalent to this one, including the property of containing
(approximately) the expected number of subsequences of each
order-type.  We provide a few applications of these results, present a construction for a family of
strongly quasirandom permutations, and prove that this construction
is essentially optimal, using a result of W. Schmidt on the
discrepancy of sequences of real numbers.
\end{abstract}

\ifx\ispreprint\undefined \indent {\it Key words:} quasirandom;
permutation; discrepancy; van der Corput
\\
\noindent {\it Correspondence Data:}

\indent Joshua N. Cooper \vspace{-.25in}

\indent UCSD Department of Mathematics \vspace{-.25in}

\indent 9500 Gilman Drive 0112 \vspace{-.25in}

\indent La Jolla, CA 92093-0112 \vspace{-.25in}

\indent jcooper@math.ucsd.edu

\pagebreak \else \fi

\section{Introduction}

In recent years, combinatorialists have been investigating several
realms of random-like -- ``quasirandom'' -- objects.  For a given
probability space $\mathcal X$, the basic idea is to choose some
collection of properties that large objects in $\mathcal X$ have
almost surely, and define a sequence $\{X_i\}_{i=1}^{\infty}
\subset \mathcal X$ to be {\it quasirandom} if $X_i$ has these
properties in the limit.  Often, this approach amounts to choosing
some random variables $\eta_j$ defined on $\mathcal X$ which tend
to their expected values almost surely as $|X| \rightarrow
\infty$, and defining $X_i$ to be quasirandom when
$(\eta_0(X_i),\eta_1(X_i),\ldots) \rightarrow
(\mathbb{E}\eta_0,\mathbb{E}\eta_1,\ldots)$ sufficiently quickly.
The resulting definitions are explored by finding many such
collections of properties and showing that quasirandomness with
respect to any one of them is equivalent to all the rest -- often
rather surprisingly, since the properties may appear completely
unrelated to one another. Quasirandom graphs, hypergraphs, set
systems, subsets of $\mathbb{Z}_n$, and tournaments have all been
examined in this way.  Quasirandom families of permutations have
been defined in [\ref{MMW1}], and Gowers [\ref{G1}] has used a careful
quantitative analysis of strongly quasirandom
(``$\alpha$-uniform'', in his terminology) subsets of
$\mathbb{Z}_n$ as an integral component of his remarkable new
proof of Szemer\'edi's Theorem.  Quasirandom objects
also have applications in algorithms as deterministic substitutes
for randomly generated objects, in addition to their purely
theoretical uses.  In fact, specific types of random-like
permutations have been used already in a number of contexts.
Lagarias [\ref{L1}] constructed random-like permutations of a
$d$-dimensional array of cells in order to solve a practical
memory-mapping problem, and Alon [\ref{A1}] used ``pseudo-random''
permutations to improve on the best known deterministic
maximum-flow algorithm of Goldberg and Tarjan.  Quasirandom sequences of reals are also fundamental to the extensively studied ``quasi-Monte Carlo'' methods of numerical analysis ([\ref{N1}]). In this paper, we
add (individual) permutations to the growing list of objects for
which a formal notion of quasirandomness has been defined.

In Section 2, we discuss the concept of $\epsilon$-balance, which
weakens the quasirandomness of Chung and Graham.  It is shown to
be equivalent to several ``types'' of quasirandomness for subsets
of $\mathbb{Z}_n$, including an infinite family of eigenvalue
bounds.  Section 3 is an excursion into the realm of subsequence
statistics of permutations, a subject that has generated a good
deal of interest recently (e.g., [\ref{AF1}] and [\ref{B1}]) --
and whose roots go back at least as far as 1935 ([\ref{ES1}]).  In
Section 4, quasirandom permutations are defined as those which map
intervals to uniformly balanced sets, and we prove that this
definition is equivalent to several other random-like conditions.  Two applications are given, including a proof that random permutations have small discrepancy.
Section 5 contains a construction for a family of strongly
quasirandom permutations that generalize the classical van der Corput sequences.  We show that this construction is
essentially optimal, using a result of Schmidt on the discrepancy
of sequences of real numbers.  Finally, Section 6 concludes with
some open problems and directions for future work.

\bigskip
\section{Balanced Sets}
\setcounter{theorem}{0}

Throughout the following, we consider permutations, i.e., elements
of $S_n$, as actions on $\mathbb{Z}_n$ as well as sequences of
numbers $(\sigma(0), \sigma(1), \ldots , \sigma(n-1))$ (``one-line
notation''). When an ordering on $\mathbb{Z}_n$ is used, we mean
the one inherited from $[0,n-1] \subset \mathbb{Z}$.  If $f_i$,
$i=1,2$, is a function from a totally ordered set $A$ to a totally
ordered set $B_i$, we say that $f_1$ and $f_2$ are
\textit{isomorphic} (and write $f_1 \sim f_2$) if, for any $a_1,
a_2 \in A$, $f_1(a_1) < f_1(a_2)$ iff $f_2(a_1) < f_2(a_2)$.  Note
that this definition still makes sense when $f_1$ and $f_2$ are
defined on different sets $A_1$ and $A_2$, so long as $|A_1| =
|A_2|$ is finite and we identify them via the unique
order-isomorphism between them. Then, if $\sigma \in S_n$ and
$\tau \in S_m$, $m \leq n$, we say that $\tau$ occurs in $\sigma$
at the set $A = \{ a_i \}_{i=1}^{m} \subset \mathbb{Z}_n$ whenever
$\sigma|_A \sim \tau$.  For each $A \subset \mathbb{Z}_n$ and
permutation $\tau$, we write ${\bf X}^{\tau}_A(\sigma)$ for the
indicator random variable of the event that $\tau$ occurs in
$\sigma$ at $A$, and we write ${\bf X}^{\tau}(\sigma)$ for the
random variable that counts the number of occurrences of $\tau$ in
$\sigma$, i.e., ${\bf X}^{\tau}(\sigma) = \sum_A {\bf
X}^{\tau}_A(\sigma)$ where $A$ ranges over all subsets of
$\mathbb{Z}_n$ of cardinality $m$.

For any subset $S \subset \mathbb{Z}_n$ (or $S \subset
\mathbb{Z}$), there is a minimal representation of $S$ as a union
of intervals.  We call these intervals the \textit{components} of
$S$ and denote the number of them by $c(S)$.  Also, we adopt the
convention that the symbols for a set and the characteristic
function of that set be the same, so, for example, $S(x) = 1$ if
$x \in S$ and $S(x) = 0$ if $x \not \in S$.  Finally, for any
function from $\mathbb{Z}_n$ to $\mathbb{C}$, we write
$\tilde{f}(k)$ for the $k^{\mbox{\scriptsize th}}$ Fourier
coefficient of $f$, defined by
$$
\tilde{f}(k) = \sum_{x \in \mathbb{Z}_n} f(x) e^{-2 \pi ikx/n}.
$$
A well known alternative definition of the Fourier coefficients of
a set $S$ is the spectrum of the circulant matrix $M_S$ whose
$(i,j)$ entry is $S(i+j)$.

\bigskip

One would expect random permutations to ``jumble'' the elements on
which it acts, i.e., there should be no correlation between
proximity in $\mathbb{Z}_n$ and proximity in the image.  We can
measure proximity by means of intervals: the elements of a small
interval are all ``close'' to one another.  Thus, if we define an
interval of $\mathbb{Z}_n$ to be the projection of any interval of
$\mathbb{Z}$, a permutation $\sigma \in S_n$ will be called
``quasirandom'' if the intersection of any interval $I$ with the
image of any other interval $J$ under $\sigma$ has cardinality
approximately $|I||J|/n$, i.e., no interval contains much more or
less of the image of any other interval than one would expect if
$\sigma$ were chosen randomly.

Thus, for any two sets $S, T \subset \mathbb{Z}_n$ we define the
{\it discrepancy} of $S$ in $T$ as
$$
D_T(S) = \left | | S \cap T | - \frac{|S||T|}{n} \right |.
$$
Note that we may apply this definition to multisets $S$ and $T$,
and that it is symmetric in its arguments.  Before proceeding, we
present a simple lemma to the effect that $D$ is subadditive:

\begin{lemma} \label{subadd} If $S = A \cup B$, $A$ and $B$ disjoint, then $D_S(T) \leq D_A(T) + D_B(T)$.  If $T = C \cup D$, $C$ and $D$ disjoint, then $D_S(T) \leq D_S(C) + D_S(D)$.  That is, $D$ is subadditive in both of its arguments.
\end{lemma}
\begin{proof}
By the triangle inequality, we have
\begin{align*}
D_T(S) & = \left | \, | \, S \cap T | - \frac{|S||T|}{n} \right | \\
& = \left | \, | \, A \cap T | + | \, B \cap T | - \frac{|A||T|}{n} - \frac{|B||T|}{n} \right | \\
& \leq \left | D_T(A) \right | + \left | D_T(B) \right |.
\end{align*}
The other statement follows by symmetry.
\end{proof}

Define $D(S)$ to be the maximum of $D_J(S)$, taken over all
intervals $J \subset \mathbb{Z}_n$, and call a set $S \subset
\mathbb{Z}_n$ $\epsilon${\it -balanced} if $D(S) < \epsilon n$.
This definition of quasirandomness is implied by that of Chung and
Graham [\ref{CG2}], according to the next proposition.

\begin{prop} \label{weak} If, for all $T \subset \mathbb{Z}_n$ and all but $\epsilon n$ of $x \in \mathbb{Z}_n$, $D_{T+x}(S) < \epsilon n$, then, for all intervals $J \subset \mathbb{Z}_n$, $D_J(S) < 2 \epsilon n$.
\end{prop}
\begin{proof}
Suppose there exists an interval $J \subset \mathbb{Z}_n$ such that
$$
\left | \, |S \cap J| - \frac{|S||J|}{n} \right | \geq 2 \epsilon
n.
$$
Then, for all $x \in \mathbb{Z}_n$,
$$
\left | \, |S \cap (J+x)| - \frac{|S||J+x|}{n} \right | \geq \left
| \, |S \cap J| - \frac{|S||J|}{n} \right | - |x|.
$$
Therefore, for each $x$ with $|x| \leq \epsilon n$, $D_{J+x}(S)
\geq \epsilon n$.  Since there are at least $\epsilon n$ such
$x$'s, setting $T = J$ contradicts the hypothesis of the
proposition.
\end{proof}

It is easy to see that the set $S = \{2x \, | \, 0 \leq x \leq n-1
\} \subset \mathbb{Z}_{2n}$ is, for any $\epsilon > 0$ and
sufficiently large $n$, $\epsilon$-bounded.  However, $S \cap
(S+t)$ does not have cardinality approximately $|S|^2/n$ for
almost all $t$, i.e., it violates ``weak translation''. Therefore,
$\epsilon$-boundedness is strictly weaker than quasirandomness in
the sense of [\ref{CG2}].

We use the convention that when ``little oh'' notation is used,
convergence in $n$ alone is intended.  (That is, the convergence
is uniform in any other quantities involved.)  The following is
the main result of this section.

\begin{theorem} \label{balance} For $r \in \mathbb{Z}_n$, we define $|r|$ to be the absolute value of the unique representative of $\, r$ from the interval $(-n/2,n/2]$.  Then, for any sequence of subsets $S \subset \mathbb{Z}_n$ and choice of $\alpha > 0$, the following are equivalent: \\ [.06in]

\noindent \begin{tabular}{p{.4in}p{4.6in}}
{\bf [B]}  & (Balance) $D(S) = o(n)$. \\
\end{tabular}

\noindent \begin{tabular}{p{.4in}p{4.6in}}
{\bf [PB]} & (Piecewise Balance) For any subset $T \subset \mathbb{Z}_n$, $D_T(S) = o(nc(T))$, where $c(T)$ denotes the number of components of $T$. \\
\end{tabular}

\noindent \begin{tabular}{p{.4in}p{4.6in}}
{\bf [MB]} & (Multiple Balance) Let $kS$ denote the multiset $\{ks|s \in S\}$.  Then, for any $k \in \mathbb{Z}_n \! \setminus \! \{0\}$, $D(kS) = o(n|k|)$. \\
\end{tabular}

\noindent \begin{tabular}{p{.4in}p{4.6in}}
{\bf [E($\frac{1}{2}$)]} & (Eigenvalue Bound $\frac{1}{2}$) For all nonzero $k \in \mathbb{Z}_n$, $\tilde{S}(k) = o(n|k|^{1/2})$. \\
\end{tabular}

\noindent \begin{tabular}{p{.4in}p{4.6in}}
{\bf [E($\alpha$)]} & (Eigenvalue Bound $\alpha$) For all nonzero $k \in \mathbb{Z}_n$, $\tilde{S}(k) = o(n|k|^{\alpha})$. \\
\end{tabular}

\noindent \begin{tabular}{p{.4in}p{4.6in}}
{\bf [S]} & (Sum) $ \sum_{r \neq 0} \left( |\tilde{S}(k)|/|k| \right )^2 = o(n^2). $ \\
\end{tabular}

\noindent \begin{tabular}{p{.4in}p{4.6in}} {\bf [T]} &
(Translation) For any interval $J$,
$$ \sum_{k \in \mathbb{Z}_n} \left ( |S \cap (J+k)| - \frac{|S||J|}{n} \right )^2 = o(n^3). $$
\end{tabular}
\end{theorem}

\ifx\ispreprint\undefined
\else
\begin{center}
\begin{figure}[ht]
\hspace{.7in} \input{figa.pstex_t} \caption{Diagram of
implications for Theorem \ref{balance}.}
\end{figure}
\end{center}
\fi

We will show that {\bf [B]} $\Rightarrow$ {\bf [PB]} $\Rightarrow$
{\bf [MB]} $\Rightarrow$ {\bf [E($\frac{1}{2}$)]} $\Rightarrow$
{\bf [E($\alpha$)]} $\Rightarrow$ {\bf [S]} $\Rightarrow$ {\bf
[T]} $\Rightarrow$ {\bf [B]}.  In each case, a statement involving
some $\epsilon$ is shown to imply the next for some $f(\epsilon)$,
where $f$ is a function which tends to zero as its argument does.
For example, Proposition \ref{division} below states that if
$D_T(S) < \epsilon n c(T)$ for all $T$, then $D_T(kS) < 2 \epsilon
n |k|$ for all k, so that $f(\epsilon)=2\epsilon$.  It appears to
be theoretically useful to track what happens to $\epsilon$ as we
pass through each implication -- see, for example, [\ref{G1}].
Thus, we include Figure 1 as an accompaniment to Theorem
\ref{balance}.  (Note that, by the proof of Proposition
\ref{MB2EA}, Figure 1 is only valid for $\epsilon < \pi / 8$,
though this is hardly a significant restriction.)  The shortcut
edge from {\bf [E($\frac{1}{2}$)]} to {\bf [S]} is given to
illustrate the (best possible) choice of $\alpha = 1/4$ in {\bf
[E($\alpha$)]}, and the circular arrow represents one complete
traversal of the cycle of implications, including the shortcut
edge.

Theorem \ref{balance} is proven in pieces, beginning with the
following proposition.

\begin{prop} {\bf [B]} $\Rightarrow$ {\bf [PB]} $\Rightarrow$ {\bf [MB]}. \label{division}
\end{prop}
\begin{proof}
Suppose that $D(S) < \epsilon n$.  Then, by Lemma \ref{subadd},
for any $T$, $D_T(S) \leq \sum D_{T_i}(S)$, where the sum is over
the components of $T$.  Thus, $D_T(S) < \epsilon n \, c(T)$, and
{\bf [B]} $\Rightarrow$ {\bf [PB]}.

Now, suppose {\bf [PB]} holds for $S$.  Note that, for a given $k
\in \mathbb{Z}_n \setminus \{0\}$ and interval $J$, the set
$J^\prime$ of elements $x \in \mathbb{Z}_n$ such that $kx \in J$
has at most $|k|$ components.  Let $J_i$ be the set of integer
points (viewed as elements of $\mathbb{Z}_n$) lying in $[a/k,b/k]
+ in/k$, so that $J^\prime = \bigcup_i J_i$.  Then the cardinality
of $J_i$ is off from $|J|/k$ by at most $1$.  By {\bf [PB]} and
the triangle inequality,
\begin{align*}
D_J(kS) &= \left | \, |k(S-t) \cap J| - \frac{|I||J|}{n} \right | \\
& \leq \left | \, \sum_i |S \cap J_i| - \sum_i \frac{|I||J_i|}{n} \right | + \frac{|I|}{n} \left | \, \sum_i |J_i| - |J| \right | \\
& < \epsilon n |k| + |k| \frac{|I|}{n} \leq 2 \epsilon n |k|.
\end{align*}
since, trivially, $\epsilon \geq n^{-1}$.
\end{proof}

Now, we wish to show that Multiple Boundedness implies the first
eigenvalue bound.  The basic idea is to imbed the elements of $S$
into the unit circle via the exponential map, and then show that a
great deal of cancellation happens because of the relatively
uniform distribution of elements of $S$.

\begin{prop} {\bf [MB]} $\Rightarrow$ {\bf [E($\frac{1}{2}$)]}. \label{MB2EA}
\end{prop}
\begin{proof}
Let $\omega = e^{2 \pi i / n}$ and $J_m^j = [\frac{nj}{m},\frac{n(j+1)}{m})$, and let $\epsilon$ be the bound on $(nk)^{-1} D(k S)$.  Recall that $\epsilon \geq n^{-1}$.  First we prove the following:
\begin{claim} Let $m$ and $j$ be positive integers with $0 \leq j < m$, and $m \geq 2$.  If we define the multiset $S_j = k S \cap J_m^j$ and let $\gamma_j = \omega^{-n (j+1/2)/m}$, then
$$
\left | \, \sum_{x \in S_j} \omega^{-x} - \frac{|S|}{m} \gamma_j
\right | < \frac{\pi |S|}{m^2} + 2 \epsilon |k|n
$$
\end{claim}

\begin{poc} We may write the left-hand side of the above expression as
\begin{align*}
\left | \, \sum_{x \in S_j} \omega^{-x} - \frac{|S|}{m} \gamma_j \right | &= \left | \, \sum_{x \in S_j} (\omega^{-x} - \frac{|S|}{m|S_j|} \gamma_j ) \right | \\
& \leq \frac{|S|}{m|S_j|} \left | \, \sum_{x \in S_j} (\omega^{-x} - \gamma_j ) \right | + \left | \, \sum_{x \in S_j} \omega^{-x} (1 - \frac{|S|}{m|S_j|}) \right | \\
& \leq \frac{|S|}{m|S_j|} \sum_{x \in S_j} \left | \, (\omega^{-x}
- \gamma_j ) \right | + \sum_{x \in S_j} \left | \, \omega^{-x} (1
- \frac{|S|}{m|S_j|}) \right |
\end{align*}
Now, for $x \in S_j$,
$$
| \, \omega^{-x} - \gamma_j | \leq | \, \omega^{-n j / m} -
\omega^{-n (j+1/2) / m} | \leq \frac{n/2}{m} \cdot \frac{2 \pi}{n}
= \frac{\pi}{m}
$$
Plugging this expression in and applying {\bf [MB]}, we have
\begin{align*}
\left | \, \sum_{x \in S_j} \omega^{-x} - \frac{|S|}{m} \gamma_j \right | & \leq \frac{|S|}{m|S_j|} \cdot |S_j| \cdot \frac{\pi}{m} + \left | |S_j| - \frac{|S|}{m} \right | \\
& \leq \frac{\pi |S|}{m^2} + \left | |kS \cap J_m^j| - \frac{|I||J_m^j|}{n} \right | + \frac{|I|}{n} \left | \frac{n}{m} - |J_m^j| \right | \\
& < \frac{\pi |S|}{m^2} + \epsilon |k| n + \frac{|I|}{n} \\
& \leq \frac{\pi |S|}{m^2} + 2 \epsilon |k| n
\end{align*}
thus, proving the claim.
\end{poc}

\medskip

If we sum over all $j \in [0, m-1)$,
\begin{align*}
\left | \, \sum_{x \in S} \omega^{-kx} \right | &= \left | \, \sum_{j=0}^{m-1} \sum_{x \in S_j} \omega^{-x} \right | \\
& \leq \left | \, \sum_{j=0}^{m-1} \frac{|S|}{m} \gamma_j \right | + \sum_{j=0}^{m-1} \left | \, \sum_{x \in S_j} \omega^{-x} - \frac{|S|}{m} \gamma_j \right | \\
& < 0 + \frac{\pi |S|}{m} + 2 \epsilon |k| nm
\end{align*}
if we assume that $m \geq 2$.  Thus, letting $m = \left \lfloor
\left ( \frac{\pi |S|}{2\epsilon |k|n} \right )^{1/2} \right
\rfloor$, we have
$$
\left | \, \sum_{x \in S} \omega^{-kx} \right | < \sqrt{18 \pi
\epsilon n|k||S|} \leq n \sqrt{18 \pi \epsilon |k|}
$$
unless $m < 2$, i.e., $ \epsilon > \frac{\pi}{8} $, which is
eventually impossible, given {\bf [MB]}.  We may therefore
conclude that $| \tilde{S}(k) | = o(n|k|^{1/2})$.
\end{proof}

A small improvement to the constant in the bound above is possible
by letting $m$ be rational, instead of integral. However, doing so
adds some complexity to the proof without making any significant
improvements.

Before we proceed with the next implication, the following lemma
will be necessary.  It implies, surprisingly, that {\bf
[E($\alpha$)]} is equivalent to {\bf [E($\beta$)]} for all
$\alpha$ and $\beta$.

\begin{lemma} \label{EA2EB} {\bf [E($\alpha$)]} implies {\bf [E($\beta$)]} for any $\alpha, \beta > 0$.
\end{lemma}
\begin{proof}
Let $M = \left \lceil \frac{\alpha}{\beta} \right \rceil$, and
assume {\bf [E($\alpha$)]}.  Then
$$
\left | \tilde{S}(k) \right |^{M} = \left | \sum_{t \in
\mathbb{Z}_n} S(t) \omega^{-kt} \right |^{M} = \left |
\sum_{t_1,\ldots,t_{M}} \left [ \prod_{j=1}^{M} S(t_j) \right ]
\omega^{-k \sum_{i=1}^M t_i} \right |
$$
Letting $u = \sum_{i=2}^M t_i$, we have
\begin{align*}
\left | \tilde{S}(k) \right |^{M} &= \left | \sum_{t_2,\ldots,t_{M}} \left [ \prod_{j=2}^{M} S(t_j) \right ] \sum_{t_1} S(t_1) \omega^{-k (t_1 + u)} \right | \\
& \leq \sum_{t_2,\ldots,t_{M}} \left [ \prod_{j=2}^{M} S(t_j) \right ] \left | \sum_{t_1} S(t_1) \omega^{-kt_1} \right | \\
&= \sum_{t_2,\ldots,t_{M}} \left [ \prod_{j=2}^{M} S(t_j) \right ] \left | \tilde{S}(k) \right | \\
&< \sum_{t_2,\ldots,t_{M}} \left [ \prod_{j=2}^{M} S(t_j) \right ] \epsilon n |k|^\alpha \\
&= |S|^{M-1} \epsilon n |k|^\alpha \leq \epsilon n^{M} |k|^\alpha.
\end{align*}
Thus, taking the $M^{\mbox{\scriptsize th}}$ root of both sides,
we have
$$
\left | \tilde{S}(k) \right | < \epsilon^{1/M} n |k|^{\alpha/M}
\leq \epsilon^{\left \lceil \alpha/\beta \right \rceil^{-1}} n
|k|^{\beta}.
$$
\end{proof}

The following corollary is actually what is needed for Theorem
\ref{balance}.

\begin{cor} {\bf [E($\frac{1}{2}$)]} $\Rightarrow$ {\bf [E($\alpha$)]}.
\end{cor}

Note that, to proceed with the next proposition, $\alpha = 1/2$
would not quite be enough -- we {\it have} to reduce it by a bit
with Proposition \ref{EA2EB}.

\begin{prop} {\bf [E($\alpha$)]} $\Rightarrow$ {\bf [S]}.
\end{prop}
\begin{proof}
By Proposition \ref{EA2EB}, we know that $|\tilde{S}(k)| <
\epsilon^{\left \lceil 4\alpha \right \rceil^{-1}} n |k|^{1/4}$
for all $k \neq 0$. Then
$$
\sum_{r \neq 0} \left( \frac{|\tilde{S}(k)|}{|k|} \right )^2 <
\sum_{k \neq 0} \left( \frac{\epsilon^{\left \lceil 4\alpha \right
\rceil^{-1}} n |k|^{1/4}}{|k|} \right )^2 \leq \epsilon^{2 \left
\lceil 4\alpha \right \rceil^{-1}} n^2 \sum_{k \neq 0} |k|^{-3/2}
< 6 \epsilon^{2 \left \lceil 4\alpha \right \rceil^{-1}} n^2.
$$
where we have used the approximation $|\zeta(s)| < (Re(s)-1)^{-1}
+ 1$ for $s$ with $Re(s) > 1$.
\end{proof}

We now write a cyclic sum in terms of Fourier coefficients. A proof of the following standard lemma is included for the sake of completeness.

\begin{lemma} \label{intb} If $J$ is an interval of $\mathbb{Z}_n$, then $\tilde{J}(k) \leq \frac{n}{2|k|}$.
\end{lemma}
\begin{proof}
We may write the magnitude of the $k^{\mbox{\scriptsize th}}$
Fourier coefficient of $J = [a+1,a+M]$ as
\begin{align*}
|\tilde{J}(k)| &= | \sum_x J(x) \omega^{-kx} | = | \sum_{x=a}^b \omega^{-kx} | = | \sum_{x=1}^M \omega^{-kx} | \\
& = \frac{|\omega^{-kM} - 1|}{|\omega^{-k} - 1|} \leq \frac{2}{4
|k| / n} = \frac{n}{2|k|}
\end{align*}
since $|e^{i \theta} - 1| \geq \frac{2 |\theta|}{\pi}$ for all
$\theta$.
\end{proof}

\begin{prop} {\bf [S]} $\Rightarrow$ {\bf [T]}.
\end{prop}
\begin{proof}
Assume that $\sum_{k \neq 0} \left( \frac{|\tilde{S}(k)|}{|k|}
\right )^2 < \epsilon n^2$.  We may write the ``translation'' sum
as
\begin{equation} \label{eqa1}
\sum_{k \in \mathbb{Z}_n} \left ( |S \cap (J+k)| -
\frac{|S||J|}{n} \right )^2 = \sum_{k \in \mathbb{Z}_n} |S \cap
(J+k)|^2 - \frac{|S|^2|J|^2}{n}
\end{equation}
Recall that $M_S$ is the $n \times n$ matrix whose $(i,j)$ entry
is $S(i+j)$.  Letting {\boldmath $v$} be the vector
$(J(0),J(1),\ldots)$, we find that $M_S \mbox{\boldmath $v$}$ is
the vector whose $k^{\mbox{\scriptsize th}}$ entry is $|I \cap
(J+k)|$.  Therefore, letting $\phi_k =
(1,\omega^{k},\omega^{2k},\ldots)$ be the $k^{\mbox{\scriptsize
th}}$ eigenvector of $M_S$,
\begin{align*}
\sum_{k \in \mathbb{Z}_n} |S \cap (J+k)|^2 &= | M_S \mbox{\boldmath $v$} |^2 = | M_S \sum_k \frac{<\mbox{\boldmath $v$},\phi_k>}{|\phi_k|^2}\,\phi_k |^2 \\
&= \sum_k | \tilde{S}(k)^2 \frac{<\mbox{\boldmath $v$},\phi_k>^2}{|\phi_k|^2} | \\
&= \sum_{k \neq 0} \left | \tilde{S}(k)
\frac{\tilde{J}(-k)}{\sqrt{n}} \right |^2 + \frac{|S|^2|J|^2}{n}
\end{align*}
Applying this equality, property {\bf [S]}, and Lemma \ref{intb}
to Equation \ref{eqa1},
\begin{align*}
\sum_{k \in \mathbb{Z}_n} \left ( |S \cap (J+k)| - \frac{|S||J|}{n} \right )^2 &= \sum_{k \neq 0} \left | \tilde{S}(k) \frac{\tilde{J}(-k)}{\sqrt{n}} \right |^2 \\
& \leq \frac{n}{4} \sum_{k \neq 0} \left | \frac{\tilde{S}(k)}{|k|} \right |^2 \\
& < \frac{\epsilon^2}{4} n^3.
\end{align*}
\end{proof}

To complete the circle of implications and finish the proof of
Theorem \ref{balance}, we show that $\epsilon$-boundedness is
implied by the ``translation'' property.

\begin{prop} {\bf [T]} $\Rightarrow$ {\bf [B]}
\end{prop}
\begin{proof}
Suppose that, for some interval $J \subset \mathbb{Z}_n$,
$$
\left | |S \cap J| - \frac{|S||J|}{n} \right | \geq 2
\epsilon^{1/3} n
$$
Then, following the line of argument given for Proposition
\ref{weak}, we may conclude that
$$
\left | |S \cap (J+k)| - \frac{|S||J|}{n} \right | \geq
\epsilon^{1/3} n
$$
whenever $|k| \leq \epsilon^{1/3} n$.  Since there are at least
$\epsilon^{1/3} n $ such $k$'s,
$$
\sum_k \left ||S \cap (J+k)| - \frac{|S||J|}{n} \right |^2 \geq
\epsilon^{1/3} n \cdot \epsilon^{2/3} n^2 = \epsilon n^3
$$
contradicting {\bf [T]}.
\end{proof}

\bigskip
\section{Statistics of Sub-Permutations}
\setcounter{theorem}{0}

Before we formally define quasirandom permutations, an excursion
into the realm of subsequence statistics is necessary. We wish to
relate ${\bf X}^{\tau}(\sigma)$, for $\tau \in S_m$, to the
quantities ${\bf X}^{\tau^\prime}(\sigma)$, with $\tau^\prime \in
S_{m+1}$, by counting occurrences of $\tau$ inside each occurrence
of $\tau^\prime$.  Define ${\bf v}_m(\sigma) \in \mathbb{Z}^{S_m}$
to be the vector whose $\tau$ component is ${\bf
X}^{\tau}(\sigma)$, and write ${\bf \tilde{v}}_m(\sigma)$ for the
vector
$$
{\bf v}_m(\sigma) - \mathbb{E}{\bf v}_m = {\bf v}_m(\sigma) -
\frac{\bf \hat{1}}{m!} \binom{n}{m}
$$
Also, let $B_m$ be the $m!$ by $(m+1)!$ matrix whose
$(\tau,\tau^\prime)$ entry is ${\bf X}^{\tau}(\tau^\prime)$ for
$\tau \in S_m$ and $\tau^\prime \in S_{m+1}$, and define $A_m$ to
be $B^*_m B_m$.

\begin{prop} For any $\sigma \in S_m$,
\begin{equation} \label{mainform}
(n-m)^2 {| {\bf \tilde{v}}_m(\sigma) |}^2 = {{\bf
\tilde{v}}_{m+1}(\sigma)}^* A_m {\bf \tilde{v}}_{m+1}(\sigma)
\end{equation}
\end{prop}
\begin{proof}
Let $\Gamma$ be the set of pairs $(U,V)$, with $U \subset V
\subset \mathbb{Z}_n$, $|U|=m$, $|V|=m+1$, and $\sigma_U \sim
\tau$.  Then, conditioning on the order-type of $U$ yields
$$
|\Gamma| = \sum_{\tau^\prime \in S_{m+1}} {\bf
X}^{\tau}(\tau^\prime) {\bf X}^{\tau^\prime}(\sigma)
$$
because each set $U$ contributes ${\bf X}^{\tau}(\sigma|_U)$ to
$|\Gamma|$.  If we instead condition on $V$ itself, then
$$
|\Gamma| = (n - m) {\bf X}^{\tau}(\sigma)
$$
because each subset $V$ is contained in exactly $n-m$ supersets
$U$.  We may therefore write
$$
\sum_{\tau^\prime \in S_{m+1}} {\bf X}^{\tau}(\tau^\prime) {\bf
X}^{\tau^\prime}(\sigma) = (n - m) {\bf X}^{\tau}(\sigma),
$$
i.e., $(n-m) {\bf v}_m(\sigma) = B_m {\bf v}_{m+1}(\sigma)$.  The
desired result then follows by linearity of expectation.
\end{proof}

Now that a numerical relationship between subsequences of length
$m+1$ and subsequences of length $m$ has been established, we need
a bound on the eigenvalues of $A_m$.

\begin{prop} \label{propsums} The following hold for all $m \geq 1$:
\begin{enumerate}
\item The column sums of $B_m$ are equal to $m+1$.
\item The row sums of $B_m$ are equal to $(m+1)^2$.
\item The row (and column) sums of $A_m$ are equal to $(m+1)^3$.
\end{enumerate}
\end{prop}
\begin{proof}
The proofs are all straightforward manipulations.

\begin{enumerate}
\item Let $b_{\tau \tau^\prime}$ denote the $(\tau,\tau^\prime)$ entry of the matrix $B_m$.  Then, denoting the set of subsets of $\mathbb{Z}_{m+1}$ of cardinality $m$ by $P_m$,
\begin{align*}
\sum_{\tau \in S_{m}} b_{\tau \tau^\prime} & = \sum_{\tau \in S_{m}} {\bf X}^{\tau}(\tau^\prime) = \sum_{\tau \in S_{m}} \sum_{A \in \mathcal{P}_m} {\bf X}^{\tau}_A(\tau^\prime) \\
& = \sum_{A \in \mathcal{P}_m} \sum_{\tau \in S_{m}} {\bf
X}^{\tau}_A(\tau^\prime) = \sum_{A \in \mathcal{P}_m} 1 = m+1.
\end{align*}
\item Note that, for a set $A \in P_m$, a permutation $\tau^\prime \in S_{m+1}$ is uniquely determined by its restriction to $A$.  Therefore, for a given $\tau \in S_m$, the number of $\tau^\prime \in S_{m+1}$ such that $\tau^\prime|_A \sim \tau$ is equal to the number of possible sets $\tau^\prime(A)$, i.e.,
$$
\sum_{\tau^\prime \in S_{m+1}} {\bf X}^{\tau}_A(\tau^\prime) =
m+1.
$$
Now, we sum over the first index:
\begin{align*}
\sum_{\tau^\prime \in S_{m+1}} b_{\tau \tau^\prime} & = \sum_{\tau^\prime \in S_{m+1}} {\bf X}^{\tau}(\tau^\prime) = \sum_{\tau^\prime \in S_{m+1}} \sum_{A \in \mathcal{P}_m} {\bf X}^{\tau}_A(\tau^\prime) \\
& = \sum_{A \in \mathcal{P}_m} \sum_{\tau^\prime \in S_{m+1}} {\bf
X}^{\tau}_A(\tau^\prime) = \sum_{A \in \mathcal{P}_m} (m+1) =
(m+1)^2.
\end{align*}
\item Let $a_{\tau \tau^\prime}$ denote the $(\tau,\tau^\prime)$ entry of the matrix $A_m$.  Since $A_m$ is symmetric, we need only show the result for column sums.
\begin{align*}
\sum_{\tau \in S_{m}} a_{\tau \tau^\prime} & = \sum_{\tau \in S_{m}} \sum_{\tau^{\prime \prime} \in S_{m+1}} b_{\tau \tau^{\prime \prime}} b_{\tau^\prime \tau^{\prime \prime}} = \sum_{\tau^{\prime \prime} \in S_{m+1}} (\sum_{\tau \in S_{m}} b_{\tau \tau^{\prime \prime}}) b_{\tau^\prime \tau^{\prime \prime}} \\
& = \sum_{\tau^{\prime \prime} \in S_{m+1}} (m+1) b_{\tau^\prime
\tau^{\prime \prime}} = (m+1)^3.
\end{align*}
where the third equality follows from part (1) and the fourth from
part (2).
\end{enumerate}
\end{proof}

\begin{cor} \label{eigcor} ${\bf \hat{1}}_m$ is an eigenvector of $A_m$ with eigenvalue $(m+1)^3$.
\end{cor}

\begin{proof}
By Proposition \ref{propsums},
\begin{align*}
A_m {\bf \hat{1}}_m = (m+1)^3 {\bf \hat{1}}_m.
\end{align*}
\end{proof}

\begin{prop} $(m+1)^3$ is the largest eigenvalue of $A_m$.
\end{prop}

\begin{proof}
By Corollary \ref{eigcor} and the Perron-Frobenius Theorem, we
need only show that $A_m$ is irreducible.  Consider the weighted
bipartite graph $G_m$ on the sets $S_m$ and $S_{m+1}$, where
$\sigma \in S_{m}$ is connected by an edge to $\tau \in S_{m+1}$
with weight ${\bf X}^\sigma (\tau)$.  (In particular, there is an
edge connecting $\sigma$ to $\tau$ iff $\sigma$ occurs in $\tau$.)
Then the adjacency matrix of $G_m$ is $B_m$, and the entries of
$A_m$ represent sums of weighted length-$2$ paths from $S_{m+1}$
to itself.  If $G_m$ is connected, then $A_m$ is irreducible.  To
establish connectivity, we show that there is a path from every
permutation $\tau \in S_{m+1}$ to the identity element of
$S_{m+1}$ in $G_m$.

\begin{claim} For each $k$, $0 \leq k \leq m$, there is a path in $G_m$ from each permutation $\tau \in S_{m+1}$ to some permutation $\tau^\prime \in S_{m+1}$ such that $\tau^\prime (i) = i$ whenever $0 \leq i \leq k$.
\end{claim}

\begin{poc}
We proceed by induction.  Suppose the claim is true for $k$, and
let $\tau$ be any element of $S_{m+1}$.  The inductive hypothesis
supplies us with a path from $\tau$ to a $\tau^\prime$ such that
$\tau^\prime (i) = i$ whenever $1 \leq i \leq k$.  Let $\sigma$ be
the unique permutation in $S_m$ such that $\sigma \sim \tau^\prime
|_{\mathbb{Z}_m \setminus \{k+1\}}$.  Then define $\tau^{\prime
\prime}$ as follows:
$$
  \tau^{\prime \prime}(i) = \left\{
               \begin{array}{ll}
               i & \text{if $0 \leq i \leq k+1$}, \\
               \sigma(i-1) + 1 & \text{if $k+1 < i \leq m$}
               \end{array}
               \right.
$$
It is easy to check that $\tau^{\prime \prime}$ is actually an
element of $S_{m+1}$.  Furthermore, $\tau^\prime|_{\mathbb{Z}_m
\setminus \{k+1\}} \sim \sigma \sim \tau^{\prime \prime}
|_{\mathbb{Z}_m \setminus \{k+1\}}$, so there is a path in $G_m$
from $\tau$ to a permutation which agrees with the identity on
$[k+1]$.
\end{poc}

Since every vertex of $S_{m}$ is connected to \emph{some} vertex
of $S_{m+1}$, and they are all connected to the identity, $G_m$ is
connected and $A_m$ is irreducible.
\end{proof}

We apply this result and the Courant-Fischer Theorem to Equation
(\ref{mainform}):
$$
(n-m)^2 {| {\bf \tilde{v}}_m(\sigma) |}^2 \leq \lambda_{max}(A_m)
| {\bf \tilde{v}}_{m+1}(\sigma) |^2 = (m+1)^3 | {\bf
\tilde{v}}_{m+1}(\sigma) |^2
$$
Thus, we have
\begin{cor} \label{down} For any $\sigma \in S_n$,
$$
| {\bf \tilde{v}}_{m}(\sigma) |^2 \leq \frac{(m+1)^3}{(n-m)^2} {|
{\bf \tilde{v}}_{m+1}(\sigma) |^2 }
$$
\end{cor}

Although Corollary \ref{down} is all we will need to use later, we
include here a short proof that $B_m$ has maximal rank.

\begin{definition}
For a permutation $\sigma \in S_{m}$, write $\sigma^\circ$ for the
permutation in $S_{m+1}$ such that $\sigma^\circ(0) = 0$ and
$\sigma^\circ(i) = \sigma(i-1) + 1$ for $1 \leq i \leq m-1$.
\end{definition}

Recall that, for two permutations $\sigma$ and $\sigma^\prime$ on
the totally ordered sets $S$ and $S^\prime$, respectively, where
$S$ and $S^\prime$ have the same cardinality, we may compare
$\sigma$ and $\sigma^\prime$ in the lexicographic order by
identifying the $i^{th}$ elements of $S$ and $S^\prime$ for each
$i$.

\begin{theorem} \label{firstoccur}
For a given $\sigma \in S_{m}$ the lexicographically least
permutation $\tau \in S_{m+1}$ such that ${\bf X}^{\sigma}(\tau)$
is nonzero is $\tau = \sigma^\circ$.
\end{theorem}
\begin{proof}
First, note that ${\bf X}^{\sigma}(\sigma^\circ) > 0$, since
$\sigma \sim \sigma^\circ |_S$, where $S = {\mathbb{Z}_{m+1}
\setminus \{0\}}$.  To show that ${\bf X}^{\sigma}(\tau) = 0$ for
all $\tau < \sigma^\circ$, we proceed by induction.  For $m=1$,
the result is trivial: $(0)$ occurs in $(01)$, and in no earlier
permutation, since $(01)$ is lexicographically first.  Now suppose
the result is true for $m-1$, but there is some $\sigma \in S_{m}$
and  $\tau \in S_{m+1}$ such that ${\bf X}^{\sigma}(\tau) \neq 0$
and $\tau < \sigma^\circ$.  We know that $\tau(0)$ must equal $0$,
since $\sigma^\circ$ does, and $\tau$ precedes $\sigma^\circ$.
Therefore, $\tau < \sigma^\circ \Rightarrow \tau |_S <
\sigma^\circ |_S$.  Since $\sigma^\circ |_S$ is isomorphic to
$\sigma$, $\tau |_S$ is lexicographically precedent to $\sigma$,
so $\sigma \not \sim \tau |_S$.  Thus, $\sigma \sim \tau |_T$ for
some $T \neq S$, i.e., a $T \in P_m$ which includes $0$.
Restricting both of these permutations to all but the first
element, we have $\sigma|_U \sim \tau|_{T \cap S} \Rightarrow {\bf
X}^{\sigma |_U}(\tau |_S) > 0$.  On the other hand, $\sigma \sim
\tau |_T \Rightarrow \sigma(0)=0$, since $0 \in T$, so the fact
that $\tau|_S < \sigma^\circ |_S \sim \sigma \sim
(\sigma|_U)^\circ$ implies that $\sigma|_U$ occurs in a
permutation lexicographically precedent to $(\sigma|_U)^\circ$,
contradicting the inductive hypothesis.
\end{proof}

\begin{cor} \label{bcor}
$\text{rank}(B_m) = m!$.
\end{cor}

\begin{proof}
If we order the columns and rows of $B_m$ by the lexicographic
order on their indices, then Theorem \ref{firstoccur} implies that
the first nonzero entry in the $i^{th}$ row occurs in the $i^{th}$
column, since $\sigma < \tau$ iff $\sigma^\circ < \tau^\circ$, and
the permutations of the form $\sigma^\circ$ precede all others.
\end{proof}

\bigskip
\section{Quasirandom Permutations}
\setcounter{theorem}{0}

In this section, we discuss several equivalent formulations of
quasirandom permutations.  The central definition is, roughly,
that a quasirandom permutation is one which sends each interval to
a highly balanced set.  Thus, we will write $D(\sigma)$ for
$\max(D_J(\sigma(I))$, where the maximum is taken over all
intervals $I$ and $J$, and a sequence of permutations ${\sigma_j}$
will be called {\it quasirandom} if $D(\sigma_j) = o(n)$.  The
following is the main result of this section.

\begin{theorem} \label{qrp} For any sequence of permutations $\sigma \in S_n$ and integer $m \geq 2$ with $n > m$, the following are equivalent: \\ [.06in]
\noindent \begin{tabular}{p{.4in}p{4.7in}}
{\bf [UB]} & (Uniform Balance) $D(\sigma) = o(n)$. \\
\end{tabular}

\noindent \begin{tabular}{p{.4in}p{4.7in}} {\bf [SP]} &
(Separability) For any intervals $I,J,K,K^\prime \subset
\mathbb{Z}_n$,
$$
\left | \sum_{x \in K \cap \sigma^{-1}(K^\prime)} I(x)
J(\sigma(x)) - \frac{1}{n} \sum_{x \in K,y \in K^\prime} I(x)J(y)
\right | = o(n)
$$ \\
\end{tabular}

\noindent \begin{tabular}{p{.4in}p{4.7in}} {\bf [mS]} &
(m-Subsequences) For any permutation $\tau \in S_m$ and intervals
$I,J \subset \mathbb{Z}_n$,
$$
{\bf X}^\tau(\sigma |_{I \cap \sigma^{-1}(J)}) = \frac{1}{m!}
\binom{|\sigma(I) \cap J|}{m} + o(n^m).
$$ \\
\end{tabular}

\noindent \begin{tabular}{p{.4in}p{4.7in}} {\bf [2S]} &
(2-Subsequences) For any intervals $I,J \subset \mathbb{Z}_n$,
$$
{\bf X}^{(01)}(\sigma |_{I \cap \sigma^{-1}(J)}) - {\bf
X}^{(10)}(\sigma |_{I \cap \sigma^{-1}(J)}) = o(n^2).
$$ \\
\end{tabular}

\end{theorem}

It follows immediately that these conditions are also equivalent
to each interpretation of the statement ``For all intervals $J
\subset \mathbb{Z}_n$, $\sigma(J)$ is $\epsilon$-balanced'' given
by the equivalences of Theorem \ref{balance}.  Thus, we have a
total of ten equivalent quasirandom properties: seven arising as
``uniformly convergent'' versions of the properties in Theorem
\ref{balance} and three new ones, which are included with uniform
balance in Figure 2.

\ifx\ispreprint\undefined
\else
\begin{center}
\begin{figure}[ht]
\hspace{1.25in} \input{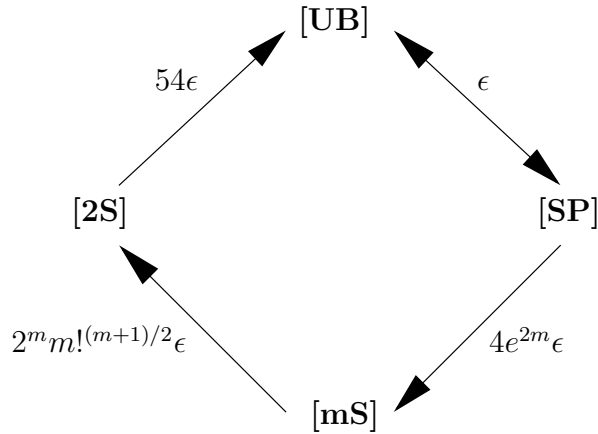} \caption{Diagram of
implications for Theorem \ref{qrp}.}
\end{figure}
\end{center}
\fi

Again, we prove the theorem piece by piece, keeping track of
$\epsilon$ as we go.  The next result states that, if uniform
balance is obeyed, then the variable $x$ and its image under
$\sigma$ are nearly independent.

\begin{prop} {\bf [UB]} $\Leftrightarrow$ {\bf [SP]}.
\end{prop}
\begin{proof}
{\bf [UB]} holds iff, for all intervals $I,J,K,K^\prime \subset
\mathbb{Z}_n$,
$$
\left | |\sigma(I \cap K) \cap (J \cap K^\prime)| - \frac{1}{n} |I
\cap K| |J \cap K^\prime| \right | < \epsilon n.
$$
But this quantity is equal to
$$
\left | \sum_{x \in K \cap \sigma^{-1}(K^\prime)} I(x)
J(\sigma(x)) - \frac{1}{n} \sum_{x \in K,y \in K^\prime} I(x) J(y)
\right |
$$
so that {\bf [UB]} is equivalent to {\bf [SP]}.
\end{proof}

Now, we show that the separability achieved in the last
proposition is sufficient to imply that subsequences happen at the
``right'' rate (i.e., what one would expect of truly random
permutations) on certain sets of indices.  A computational lemma
will greatly simplify the proof.

\begin{lemma} \label{prodest} If, for each $j$ with $1 \leq j \leq k$, $|a_j|<n^{c_j}$, $\epsilon_j > n^{-1}$, and
$$
|x_j - a_j| < \epsilon_j n^{c_j-1}
$$
then
$$
\left | \prod_{j=1}^k x_j - \prod_{j=1}^k a_j \right | < 3^{k-1}
n^{\sum_{j=1}^k c_j - k} \prod_{j=1}^k \epsilon_j.
$$
\end{lemma}
\begin{proof}
We show the result for $k=2$, and the general case follows by a
simple induction.  Thus,
\begin{align*}
\left | x_1 x_2 - a_1 a_2 \right | & \leq \left | x_1 - a_1 \right | \cdot \left | x_2 - a_2 \right | + a_1  \left | x_2 - a_2 \right | + a_2 \left | x_1 - a_1 \right | \\
&< \epsilon_1 \epsilon_2 n^{c_1 + c_2 - 2} + \epsilon_1 n^{c_1 + c_2 - 1} + \epsilon_2 n^{c_1 + c_2 - 1} \\
&< 3 \epsilon_1 \epsilon_2 n^{c_1 + c_2 - 2}.
\end{align*}
\end{proof}

\begin{prop} {\bf [SP]} $\Rightarrow$ {\bf [mS]}.
\end{prop}
\begin{proof}
Let $I,J$ be intervals, and let $K = {I \cap \sigma^{-1}(J)}$.
Note that we may write the number of ``occurrences'' of $\tau \in
S_m$ in $\sigma |_K$ as
$$
{\bf X}^\tau(\sigma |_K) = \sum_{x_1,\cdots,x_m \in K}
\prod_{i=0}^m \left ( \chi(x_i < x_{i+1})
\chi(\sigma(x_{\tau^{-1}(j)}) < \sigma(x_{\tau^{-1}(j+1)})) \right
)
$$
In the interest of notational compactness, we will denote
$\chi(x<y)$ by $\cl{x}{y}$, and define $\cl{x}{y} = 1$ if either
$x$ or $y$ is undefined.  Furthermore, for any subset $A \subset
[m]$, we will denote the following expression
$$
\sum_{\{x_i\}_{i \not \in A} \subset K} \, \sum_{\{x_k\}_{k \in A}
\subset I} \sum_{\{x^{\prime}_k\}_{k \in A} \subset J} \left (
\prod_{j=0}^m \cl{x_j}{x_{j+1}}
\cl{x^A_{\tau^{-1}(j)}}{x^A_{\tau^{-1}(j+1)}} \right )
$$
by $\Sigma(A)$, where $x^A_j$ means $\sigma(x_j)$ for $j \not \in
A$, and $x^{\prime}_k$ for $j \in A$.  Thus, ${\bf X}^\tau(\sigma
|_K) = \Sigma(\emptyset)$.  The proof will now proceed by
induction on the subsets of $[m]$, ordered by inclusion.

Suppose $A \subset B \subset [m]$, with $B \setminus A = \{s\}$,
and assume that
\begin{equation} \label{indstep}
\left | {\bf X}^\tau(\sigma |_K) - n^{-|A|} \Sigma(A) \right | <
|A| \epsilon n^m
\end{equation}
By {\bf [SP]}, we know that, for any $a,b,c,d \in \mathbb{Z}_n$,
the quantity
$$
\left | \sum_{x_s \in K} \cl{a}{x_s} \cl{x_s}{b}
\cl{c}{\sigma(x_s)} \cl{\sigma(x_s)}{d} - \frac{1}{n} \sum_{x_s
\in I, x_s^\prime \in J} \cl{a}{x_s} \cl{x_s}{b}
\cl{c}{x^\prime_s} \cl{x^\prime_s}{d} \right |
$$
is bounded above by $\epsilon n$.  Then, substituting $a=x_{s-1}$,
$b=x_{s+1}$, $c=x^A_{\tau^{-1}(\tau(s)-1)}$, and
$d=x^A_{\tau^{-1}(\tau(s)+1)}$ to account for all the terms
containing $x_s$ in the product portion of the expression
$\Sigma(A)$, we have (after a very messy but otherwise
straightforward calculation),
$$
\left | \Sigma(A) - n^{-1} \Sigma(B) \right | < \epsilon n
|K|^{m-|B|} |I|^{|A|} |J|^{|A|} \leq \epsilon n^{m + |A|}
$$
Applying this to the inductive hypothesis with the aid of the
triangle inequality yields
\begin{align*}
\left | {\bf X}^\tau(\sigma |_K) - n^{-|B|} \Sigma(B) \right | & \leq \left | {\bf X}^\tau(\sigma |_K) - n^{-|A|} \Sigma(A) \right | + n^{-|A|} \left | \Sigma(A) - n^{-1} \Sigma(B) \right | \\
& < |A| \epsilon n^m + n^{-|A|} \epsilon n^{m+|A|} = |B| \epsilon
n^m.
\end{align*}
Therefore, (\ref{indstep}) is true for all $A \subset [m]$.  In
particular, it is true for $A = [m]$, so that
$$
\left | {\bf X}^\tau(\sigma |_K) - n^{-m} \Sigma([m]) \right | < m
\epsilon n^m
$$
Since we have
\begin{align}
\nonumber \left | \sum_{\{x_j\} \subset I} \prod_{j=0}^m \cl{x_j}{x_{j+1}} - \frac{|I|^m}{m!} \right | & = \left | \binom{|I|}{m} - \frac{|I|^m}{m!} \right | \leq \frac{(|I|+m)^m - |I|^m}{m!} \\
\nonumber & = \frac{1}{m!} \sum_{k=1}^{m} \binom{m}{k} |I|^{m-k} m^k \leq \frac{|I|^{m-1}}{m!} \sum_{k=0}^{m} \binom{m}{k} m^k \\
& = \frac{(1+m)^m}{m!} |I|^{m-1} \label{approxbin}
\end{align}
and also
$$
\left | \sum_{\{x_j\} \subset J} \prod_{j=0}^m \cl{x_j}{x_{j+1}} -
\frac{|J|^m}{m!} \right | \leq \frac{(1+m)^m}{m!} |J|^{m-1}
$$
we may conclude that
$$
\left | \Sigma([m]) - \frac{|I|^m|J|^m}{m!^2} \right | < 3
n^{2m-1} \frac{(1+m)^{2m}}{m!^2}
$$
by Lemma \ref{prodest}.  Thus,
\begin{align*}
\left | {\bf X}^\tau(\sigma |_K) - \frac{1}{m!^2} \left ( \frac{|I||J|}{n} \right )^m \right | & \leq \left | {\bf X}^\tau(\sigma |_K) - \frac{1}{n^m} \Sigma([m]) \right | \\
& + \frac{1}{n^m} \left | \Sigma([m]) - \frac{|I|^m|J|^m}{m!^2} \right | \\
& < m \epsilon n^m + 3 n^{m-1} \frac{(1+m)^{2m}}{m!^2} \\
& < 4 \epsilon n^m \frac{(1+m)^{2m}}{m!^2}
\end{align*}
But, by {\bf [UB]} (which is equivalent to {\bf [SP]}), and Lemma
\ref{prodest}
$$
\left | \left ( \frac{|I||J|}{n} \right )^m - |\sigma(I) \cap J|^m \right | < 3^{m-1} \epsilon^m n^m \\
$$
Since $m \geq 2$, this gives
\begin{align*}
\left | {\bf X}^\tau(\sigma |_K) - \frac{1}{m!^2} |\sigma(I) \cap J|^m \right | &< n^m \frac{1}{m!^2} \left (3^{m-1} \epsilon^m + 4 \epsilon (1+m)^{2m} \right ) \\
& < n^m \frac{1}{m!^2} \left ((1+m)^{2m} \epsilon + 4 \epsilon (1+m)^{2m} \right) \\
& = \frac{5 \epsilon (1+m)^{2m}}{m!^2} n^m
\end{align*}
Finally, the fact that $\epsilon \geq n^{-1}$ implies, as in
(\ref{approxbin}),
$$
\left | \binom{|\sigma(I) \cap J|}{m} - \frac{|\sigma(I) \cap
J|^m}{m!} \right | \leq \epsilon \frac{(1+m)^m}{m!} n^m
$$
so we may conclude
$$
\left | {\bf X}^\tau(\sigma |_K) - \frac{1}{m!} \binom{|\sigma(I)
\cap J|}{m} \right | < \frac{6 \epsilon (1+m)^{2m}}{m!^2} n^m < 4
e^{2m} \epsilon n^m
$$
where we have used the Stirling approximation $m! > \sqrt{2 \pi m}
\left (\frac{m}{e} \right)^m$.
\end{proof}

Now, we use the results of the previous section to show that {\bf
[mS]} implies {\bf [2S]}.

\begin{prop} {\bf [mS]} $\Rightarrow$ {\bf [2S]}.
\end{prop}
\begin{proof}
Let $K=\sigma(I) \cap J$ for some intervals $I,J \in
\mathbb{Z}_n$.  We may assume that $n \geq 2m$, so that $n-k>n/2$
for all $k<m$.  Therefore, by Corollary \ref{down}, we may write
$$
| {\bf \tilde{v}}_{m-1}(\sigma|_K) |^2 \leq \frac{m^3}{(n-m+1)^2}
{| {\bf \tilde{v}}_{m}(\sigma|_K) | } < \frac{4 m^3}{n^2} m!
\max_{\tau \in S_m} \left | X^\tau(\sigma|_K) - \binom{|K|}{m}
\right |^2
$$
Iterating this process $m-2$ times, we find
\begin{align*}
| {\bf \tilde{v}}_{m-1}(\sigma) |^2 &< \frac{4^{m-5} m!^{m+1}}{n^{2m-4}} \max_{\tau \in S_m} \left | X^\tau(\sigma|_K) - \binom{|K|}{m} \right | \\
& < \frac{4^{m-2} m!^{m+1}}{n^{2m-4}} \epsilon^2 n^{2m} = 4^{m-5}
m!^{m+1} \epsilon^2 n^4.
\end{align*}
Let the quantity $d$ be defined by
$$
d = \left | X^{(01)}(\sigma|_K) - \binom{|K|}{2} \right | = \left
| X^{(10)}(\sigma|_K) - \binom{|K|}{2} \right |.
$$
Then $| {\bf \tilde{v}}_{2}(\sigma) |^2 = 2d^2$, so $d < 2^{m -
11/2} m!^{(m+1)/2} \epsilon n^2$, and
$$
\left | X^{(01)}(\sigma|_K) - X^{(10)}(\sigma|_K) \right | < 2d <
2^m m!^{(m+1)/2} \epsilon n^2
$$
which implies {\bf [2S]}.
\end{proof}

In what follows, we denote the complement of a set $S \in
\mathbb{Z}_n$ by $\bar{S}$, and we denote by $S^*$ its projection
onto $[0,n-1]$.  Also, call an interval $I \subset \mathbb{Z}_n$
``contiguous'' if $I^*$ is an interval, ``terminal'' if $\bar{I}$
is contiguous, ``initial'' if it is terminal and contains $0$, and
``final'' if it is terminal and contains $n-1$.

\begin{prop} {\bf [2S]} $\Rightarrow$ {\bf [UB]}.
\end{prop}
\begin{proof}
Suppose $\sigma$ satisfies {\bf [2S]} but not {\bf [UB]}. We claim
that, for infinitely many $n$ and some $\epsilon > 0$, there are
intervals $I, J \subset \mathbb{Z}_n$ with $I$ and $J$ initial,
and $D_J(\sigma(I))$ at least $27\epsilon n/2$.  Since {\bf [UB]}
is not true for $\sigma$, we may choose $\epsilon$ so that there
are proper subintervals $I,J \subset \mathbb{Z}_n$ with
$D_{J}(\sigma(I)) \geq 54 \epsilon n$.  Suppose $J$ is not
contiguous.  Then $\bar{J}$ is contiguous, and, since
$$
\left ( \left | \sigma(I) \cap J  \, \right | - \frac{|I||J|}{n}
\right ) + \left ( \left | \sigma(I) \cap \bar{J} \, \right | -
\frac{|I||\bar{J}|}{n} \right ) = \left | \sigma(I) \right | -
\frac{|I|n}{n} = 0
$$
we may replace $J$ with $\bar{J}$ and retain the property that
$D_{J}(\sigma(I)) \geq 54 \epsilon n$.  Now, suppose $J$ is not
terminal.  Let $J^\prime$ be a component of $(\bar{J})^*$.
$J^\prime$ is terminal because $J$ is contiguous, and, since we
have
$$
\left | \sigma(I) \cap (J \cup J^\prime)  \, \right | -
\frac{|I||J \cup J^\prime |}{n} = \left ( \left | \sigma(I) \cap J
\, \right | - \frac{|I||J|}{n} \right ) + \left ( \left |
\sigma(I) \cap J^\prime \, \right | - \frac{|I||J^\prime|}{n}
\right )
$$
either $D_{J \cup J^\prime}(\sigma(I)) \geq 27 \epsilon n$ or
$D_{J^\prime}(\sigma(I)) \geq 27 \epsilon n$.  Thus, we may assume
that $J$ is terminal (since $J^\prime$ and $J \cup J^\prime$ are),
and $D_{J}(\sigma(I)) \geq 27 \epsilon n$.  If $J$ is final,
taking its complement makes it initial without disturbing the
discrepancy.  Apply the same process to $I$ to ensure that it is
initial, with the penalty that now
\begin{equation} \label{tst}
D_{J}(\sigma(I)) \geq 27 \epsilon n/2.
\end{equation}

For ease of notation, we will let
\begin{eqnarray*}
A & = I \cap \sigma^{-1}(J) & \hspace{.3in} a=|A| \\
B & = I \cap \sigma^{-1}(\bar{J}) & \hspace{.3in} b=|B| \\
C & = \bar{I} \cap \sigma^{-1}(J) & \hspace{.3in} c=|C| \\
D & = \bar{I} \cap \sigma^{-1}(\bar{J}) & \hspace{.3in} d=|D|
\end{eqnarray*}
For subsets $S,T \subset \mathbb{Z}_n$, let $\partial_\sigma
(S,T)$ denote the number of pairs $(x,y) \in S \times T$ such that
$x < y$ and $\sigma(x) < \sigma(y)$.  Then
\begin{align*}
{\bf X}^{(01)}(\sigma) & = {\bf X}^{(01)}(\sigma|_I) + {\bf X}^{(01)}(\sigma|_{\bar{I}}) + \partial_\sigma(I, \bar{I}) \\
& = {\bf X}^{(01)}(\sigma|_I) + {\bf X}^{(01)}(\sigma|_{\bar{I}}) + \partial_\sigma(A, C) \\
& \, \, \, \, \, + \partial_\sigma(A, D) + \partial_\sigma(B, C) +
\partial_\sigma(B, D)
\end{align*}
Now, $\partial_\sigma(B, C) = 0$ and $\partial_\sigma(A, D) = ad$,
since every element of $J$ is less than every element of
$\bar{J}$, and every element of $I$ is less than every element of
$\bar{I}$.  Also,
\begin{align*}
\partial_\sigma(A, C) & = {\bf X}^{(01)}(\sigma|_{\sigma^{-1}(J)}) - {\bf X}^{(01)}(\sigma|_A) - {\bf X}^{(01)}(\sigma|_C) \\
\partial_\sigma(B, D) & = {\bf X}^{(01)}(\sigma|_{\sigma^{-1}(\bar{J})}) - {\bf X}^{(01)}(\sigma|_B) - {\bf X}^{(01)}(\sigma|_D)
\end{align*}
Thus, we have
\begin{align*}
{\bf X}^{(01)}(\sigma) & = {\bf X}^{(01)}(\sigma|_I) + {\bf X}^{(01)}(\sigma|_{\bar{I}}) + {\bf X}^{(01)}(\sigma|_{\sigma^{-1}(J)}) - {\bf X}^{(01)}(\sigma|_A) \\
& - {\bf X}^{(01)}(\sigma|_C) + ad + {\bf
X}^{(01)}(\sigma|_{\sigma^{-1}(\bar{J})}) - {\bf
X}^{(01)}(\sigma|_B) - {\bf X}^{(01)}(\sigma|_D)
\end{align*}
By {\bf [2S]}, for sufficiently large $n$, each term ${\bf
X}^{(01)}(\sigma|_S)$ can be approximated by $\binom{|S|}{2}/2$ to
within $\epsilon n^2/2$, and therefore by $|S|^2/4$ to within $3
\epsilon n^2/4$, since
$$
\left | \binom{|S|}{2} - \frac{|S|^2}{2} \right | = \frac{|S|}{2}
< \frac{\epsilon}{2} n^2.
$$
Therefore, rewriting and multiplying by $4$, we have that
$$
\left | n^2 - (a+b)^2 - (c+d)^2 - (a+c)^2 - (a+d)^2 + a^2 + b^2 +
c^2 + d^2 - 4ad \right |
$$
is bounded above by $27 \epsilon n^2$.  Since $n = a+b+c+d$, we
may simplify down to
\begin{equation} \label{eq1}
| bc - ad | < \frac{27 \epsilon}{2} n^2
\end{equation}
Let $\delta n = \left | I \cap \sigma^{-1}(J) \right | -
|I||J|/n$.  Then, by (\ref{tst}),
\begin{align*}
| bc - ad | &= \left | \left ( \frac{|I||\bar{J}|}{n} - \delta n \right ) \left ( \frac{|\bar{I}||J|}{n} - \delta n \right ) - \left ( \frac{|I||J|}{n} + \delta n \right ) \left ( \frac{|\bar{I}||\bar{J}|}{n} + \delta n \right ) \right | \\
&= |\delta| (|I||J| + |I||\bar{J}| + |\bar{I}||J| + |\bar{I}||\bar{J}|) \\
&= \frac{D_J(\sigma(I))}{n} \cdot (|I|+|\bar{I}|)(|J|+|\bar{J}|)
\geq \frac{27 \epsilon}{2} n^2
\end{align*}
contradicting (\ref{eq1}).
\end{proof}

We present two simple applications of these results.  The following observation has some relevance to the investigations
of [\ref{AF1}] and [\ref{B1}].

\begin{prop} If a permutation $\sigma \in S_n$ excludes $\tau \in
S_m$ (in the sense that ${\bf X}^\tau(\sigma) = 0$), then
$$
D(\sigma) \geq \frac{n\binom{n}{m}}{4 \, e^{2m} \,m!\, n^m}.
$$
\end{prop}
\begin{proof} Let $\epsilon = D(\sigma)/n$.  We show that, if
$$
\epsilon \leq \frac{\binom{n}{m}}{4 \, e^{2m} \, m! \, n^m}
$$
then there is at least one copy of every element $\tau$ of $S_m$
in $\sigma$.  According to the implication from {\bf [UB]} to {\bf
[mS]}, if $D(\sigma) < \epsilon n$, then
$$
\left | \textbf{X}^\tau(\sigma) - \frac{1}{m!} \binom{n}{m} \right
| < 4 e^{2m} \epsilon n^m \leq \frac{1}{m!} \binom{n}{m},
$$
so that $\textbf{X}^\tau(\sigma) > 0$.
\end{proof}

\begin{cor} There is a constant $c > 0$ so that, if $n \geq 2m$ and
$\sigma \in S_n$ excludes $\tau \in S_m$, then
$$
\frac{D(\sigma)}{n} > \left(\frac{c}{m^2} \right)^m.
$$
\end{cor}

We can also use Theorem \ref{qrp} to calculate the discrepancy of a random permutation.

\begin{theorem} If a permutation $\sigma$ is chosen randomly and uniformly from $S_n$, then $D(\sigma) = O(\sqrt{n \log n})$ almost surely.
\end{theorem}
\begin{proof} We use $c_i$, $i \in \mathbb{N}$ to denote absolute constants throughout.  Let the random variable $\xi_n$ be the number of inversions in a randomly and uniformly chosen element of $S_n$.  Define $\eta_n$ to be the normalized random variable given by
$$
\eta_n = \frac{\xi_n - \mbox{\textbf{E}}(\xi_n)}{(\mbox{Var}\, \xi_n)^{1/2}} = \frac{\xi_n - {\frac{1}{2} \binom{n}{2}}}{(\mbox{Var}\, \xi_n)^{1/2}},
$$
and let $u_n$ denote the distribution function with an atom of mass $1/n$ at each of $0, \ldots , n-1$.  It is well known ([\ref{St1}],[\ref{St2}]) that the generating function $g_n(q)$ for the number of permutations with a given number of inversions is the $q$-factorial $[n]!$, and that its coefficients are symmetric and unimodal.  Therefore, $g_n = g_{n-1} \cdot (q^{n-1} + \ldots + 1)$, so that $\xi_n = \xi_{n-1} \ast u_n$ (the convolution product).  The unimodality of the coefficients of $g_{n-1}$ implies the concavity of the cumulative distribution function $F_n(x)$ of $\xi_n$ on the interval $[\lceil \binom{n}{2}/2 \rceil,\infty]$.  Thus, if $x \geq \binom{n-1}{2}/2+n-1$, then $F_n(x) = (F_{n-1} \ast u_n)(x) \leq F_{n-1}(x)$.  (The convolution is a finite sum, so the implicit change of order of summation is legitimate.)  Similarly, $F_n(x) = (F_{n-1} \ast u_n)(x) \geq F_{n-1}(x)$ if $x \leq \binom{n-1}{2}/2$.  We have, then,
\begin{equation} \label{monotone}
\mbox{\textbf{Pr}}( \left|\xi_n - \mbox{\textbf{E}}(\xi_n)\right| > \lambda) \geq \mbox{\textbf{Pr}}( \left|\xi_{m} - \mbox{\textbf{E}}(\xi_m)\right| > \lambda )
\end{equation}
whenever $\lambda \geq (n-1)/2$ and $m \leq n$.

It is a theorem of Sachkov [\ref{Sa1}] that the cumulative distribution function of $\eta_n$ converges to $\Phi(0,1)$ (the c.d.f. of the standard normal distribution), and that $\sigma_n^2 = \mbox{Var}\, \xi_n = n^3/36 + O(n^2)$.  In particular, the moment generating function $M(t,n)=\mbox{\textbf{E}}(e^{t \eta_n})$ of $\eta_n$ is given by
$$
\log M(t,n) = \frac{t^2}{2} + \sum_{k=2}^\infty \frac{B_{2k} \, t^{2k}}{2k \, \sigma_n^{2k} \, (2k)!} \sum_{j=1}^n (j^{2k}-1)
$$
where $B_{2k}$ is the $2k^{\mbox{th}}$ Bernoulli number.  Let $f_n(t) = \log M(t,n) - t^2/2$.  Then, using the approximation $|B_{2k}| < 4 (2k)!/(2 \pi)^{2k}$, we have
\begin{equation} \label{momest}
|f_n(t)| \leq c_1 \sum_{k=2}^\infty \frac{t^{2k}}{(2 \pi)^{2k} \, 2k \, n^{3k}} \cdot n^{2k+1} \leq c_2 n \sum_{k=2}^\infty \left ( \frac{t}{2 \pi n^{1/2}} \right )^{2k} \leq \frac{c_3 t^4}{4 \pi^2 n - t^2},
\end{equation}
so long as $t < n^{1/2}$.

By Theorem \ref{qrp}, there exists an $\alpha>0$ so that
$$
\mbox{\textbf{Pr}}(D(\sigma) > \lambda \sqrt{n \log n}) \leq \mbox{\textbf{Pr}}(\max_{I,J} \left|{\bf X}^{(10)}(\sigma|_{I \cap \sigma^{-1}(J)}) - \frac{1}{2} \binom{|I \cap \sigma^{-1}(J)|}{2}\right| > c_4 \lambda n^{3/2} \sqrt{\log n}).
$$
By (\ref{monotone}), we may write
\begin{align*}
\mbox{\textbf{Pr}}(D(\sigma) > \lambda \sqrt{n \log n}) & \leq \sum_{I,J} \mbox{\textbf{Pr}}(\left|\xi_{|\sigma(I) \cap J|} - \mbox{\textbf{E}} \xi_{|\sigma(I) \cap J|}\right| > c_4 \lambda n^{3/2} \sqrt{\log n}) \\
& \leq n^4 \mbox{\textbf{Pr}}(\left|\xi_n - \mbox{\textbf{E}} \xi_n \right| > c_4 \lambda n^{3/2} \sqrt{\log n})
\end{align*}
so long as $n$ is sufficiently large.  Furthermore, by Markov's inequality and the estimate on $\sigma_n$, for fixed $\lambda > 0$,
\begin{align*}
\mbox{\textbf{Pr}}(\left|\xi_n - \mbox{\textbf{E}} \xi_n \right| > c_4 \lambda n^{3/2} \sqrt{\log n}) & \leq \mbox{\textbf{Pr}}(\left|\eta_n\right| > c_5 \lambda \sqrt{\log n}) \\
&= 2 \mbox{\textbf{Pr}}(e^{\eta_n} > e^{c_5 \lambda \sqrt{\log n}} ) \\
&\leq 2 \mbox{\textbf{E}}(e^{t \eta_n}) e^{-t c_5 \lambda \sqrt{\log n}}.
\end{align*}
Setting $t = c_5 \lambda \sqrt{\log n}$ and applying the bound (\ref{momest}), we have
$$
\mbox{\textbf{Pr}}(\left|\xi_n - \mbox{\textbf{E}} \xi_n \right| > c_4 \lambda n^{3/2} \sqrt{\log n}) \leq e^{- c_6 \lambda^2 \log n + c_7 \log^2 n/n} \leq c_8 e^{- c_6 \lambda^2 \log n}.
$$
Therefore,
$$
\mbox{\textbf{Pr}}(D(\sigma) > \lambda \sqrt{n \log n}) \leq c_8 n^{4 - c_6 \lambda^2},
$$
which tends to zero if we choose $\lambda > 2 c_6^{-1/2}$.
\end{proof}

\bigskip
\section{Constructions}
\setcounter{theorem}{0}

In this section, we present a construction for a large class of permutations which are highly quasirandom.  We will assume throughout that $\sigma \in S_n$ and $\tau \in S_m$, unless indicated otherwise.

\begin{definition} For permutations $\sigma \in S_n$ and $\tau \in S_m$, considered as actions on $\mathbb{Z}_n$ and $\mathbb{Z}_m$, respectively, define $\sigma \otimes \tau \in S_{nm}$ by $(\sigma \otimes \tau)(x) = \tau(\lfloor \frac{x}{n} \rfloor) + m \sigma(x \! \! \mod n)$.  We will also denote the $k^{th}$ product of $\sigma$ with itself as $\sigma^{(k)}$.
\end{definition}

A special case of this product appears in [\ref{D2}], where the authors define a sequence of permutations lacking ``monotone $3$-term arithmetic progressions'' by taking iterated products of the elements of $S_2$.

Note that $\sigma \otimes \tau$ has the property that $(\sigma
\otimes \tau)([0,n-1])$ is the set of all elements of
$\mathbb{Z}_{nm}$ congruent to $0$ mod $m$ (i.e., $m \cdot
[0,n-1]$), a set which necessarily lacks the ``weak translation''
property of quasirandom sets.  Thus, a sequence $\{
\sigma_1,\sigma_1 \otimes \sigma_2,\sigma_1 \otimes \sigma_2
\otimes \sigma_3,\ldots \}$ sends intervals to sets which are not
quasirandom in the sense of [\ref{CG2}].  Nonetheless, we will
prove shortly that it does satisfy ${\bf UB}$.  First, we
offer a justification for the lack of parentheses in the
expression for this sequence.

\begin{prop} $\otimes$ is associative.
\end{prop}
\begin{proof}
Suppose $\sigma_i \in S_{n_i}$ for $i=1,2,3$.  Then, applying the
definition of $\otimes$ twice,
\begin{align*}
\left [(\sigma_1 \otimes \sigma_2) \otimes \sigma_3 \right] (x) & = \sigma_3 \left( \left \lfloor \frac{x}{n_1n_2} \right \rfloor \right ) + n_3 \sigma_2 \left ( \left \lfloor \frac{x \mbox{ mod } n_1n_2}{n_1} \right \rfloor \right ) \\
& \, \, \, \, \, + n_2n_3 \sigma_1(x \mbox{ mod } n_1)
\end{align*}
and
\begin{align*}
\left [\sigma_1 \otimes (\sigma_2 \otimes \sigma_3) \right ](x) & = \sigma_3 \left( \left \lfloor \frac{\left \lfloor x/n_1 \right \rfloor}{n_2} \right \rfloor \right ) + n_3 \sigma_2 \left ( \left \lfloor \frac{x}{n_1} \right \rfloor \mbox{ mod } n_2 \right ) \\
& \, \, \, \, \, + n_2n_3 \sigma_1(x \mbox{ mod } n_1)
\end{align*}
Note that every element $x$ of $\mathbb{Z}_{n_1n_2n_3}$ can be
represented uniquely as $x = an_1n_2 + bn_1 + c$, with $0 \leq a <
n_3$, $0 \leq b < n_2$, and $0 \leq c < n_1$.  Using this
notation, we find that
\begin{align*}
\left [(\sigma_1 \otimes \sigma_2) \otimes \sigma_3 \right] (x) & = \sigma_3(a) + n_3 \sigma_2 \left ( \left \lfloor \frac{bn_1 + c}{n_1} \right \rfloor \right ) + n_2n_3 \sigma_1(c) \\
& = \sigma_3(a) + n_3 \sigma_2(b) + n_2n_3 \sigma_1(c)
\end{align*}
and
\begin{align*}
\left [\sigma_1 \otimes (\sigma_2 \otimes \sigma_3) \right ](x) = \,\, & \sigma_3 \left( \left \lfloor \frac{an_2 + b}{n_2} \right \rfloor \right ) + n_3 \sigma_2 \left ( (a n_2 + b) \mbox{ mod } n_2 \right ) \\
& \,\, + n_2n_3 \sigma_1(c) \\
= & \,\, \sigma_3(a) + n_3 \sigma_2(b) + n_2n_3 \sigma_1(c).
\end{align*}
\end{proof}

Define $d(\sigma)$ by
$$
d(\sigma) = \max_{I,J} D_J(\sigma(I))
$$
where $I$ is allowed to vary over all possible intervals, but $J$
is restricted to initial intervals.  We denote the analogue for
final intervals by $d^\prime$.  Then we have the following result:

\begin{prop} $d(\sigma \otimes \tau) \leq m - 1 + d(\sigma)$.
\end{prop}
\begin{proof}
Let the interval $I_k=[kn,(k+1)n-1] \subset \mathbb{Z}_{nm}$.
Then, any initial interval $S$ of $\mathbb{Z}_{nm}$ can, for some
$l<m$, be written
$$
S = \bigcup_{k=0}^l I_k \cup S_0
$$
where $S_0$ is an initial segment of $I_{l+1}$.  For any interval
$J \subset \mathbb{Z}_{nm}$, then, we may write
$$
D_J((\sigma \otimes \tau)(S)) \leq \sum_{k=0}^l D_J((\sigma
\otimes \tau)(I_k)) + D_J((\sigma \otimes \tau)(S_0))
$$
by Lemma \ref{subadd}.  First, we estimate $D_J(\sigma(I_k))$.
\begin{align*}
D_J((\sigma \otimes \tau)(I_k)) &= \left | \, | \, (\sigma \otimes \tau)(I_k) \cap J | - \frac{|(\sigma \otimes \tau)(I_k)||J|}{nm} \right | \\
&= \left | \, | \, (m[0,n-1]+k) \cap J | - \frac{n|J|}{nm} \right | \\
&\leq \left | \, \frac{|J| + m - 1}{m} - \frac{|J|}{m} \right | =
\frac{m-1}{m}
\end{align*}
Let $J_0 \subset \mathbb{Z}_n$ denote the set $\left \{ \lfloor
\frac{x}{m} \rfloor | \, x \in J \right \}$, and let $S_1 \subset
\mathbb{Z}_n$ be the set $S_0$ reduced mod $n$.  Then,
\begin{align*}
D_J((\sigma \otimes \tau)(S_0)) &= \left | \, | \, (\sigma \otimes \tau)(S_0) \cap J | - \frac{|(\sigma \otimes \tau)(S_0)||J|}{nm} \right | \\
&= \left | \, | \, \sigma(S_1) \cap J_0 | - \frac{|S_1||J|}{nm} \right | \\
&= \left | \, | \, \sigma(S_1) \cap J_0 | - \frac{|S_1|m|J_0|}{nm} + \frac{|S_1|}{nm} (m|J_0|-|J|) \right | \\
& \leq D_{J_0}(\sigma(S_1)) + \left ( \, \frac{n}{nm} (m-1) \right
) \leq d(\sigma) + \frac{m-1}{m}
\end{align*}
Thus,
\begin{align*}
d(\sigma \otimes \tau) & \leq (m-1) \frac{m-1}{m} + \frac{m-1}{m} + d(\sigma) \\
& = m - 1 + d(\sigma).
\end{align*}
\end{proof}

An identical result holds for $d^\prime$, by symmetry.  We use
this in the next proposition, which allows us to bound
discrepancies recursively.

\begin{prop} \label{inh} $D(\sigma \otimes \tau) \leq m-1 + d(\sigma) + d^\prime(\sigma)$.
\end{prop}
\begin{proof}
Note that every interval $I$ of $\mathbb{Z}_{nm}$ is of the form
$$
S = \bigcup_{k \in [l,L]} I_k \cup S_0 \cup S^\prime_0
$$
where $[l,L]$ is an interval of $\mathbb{Z}_m$ of length no more
than $m-2$, $S_0$ is an initial segment of $I_{L+1}$, and
$S^\prime_0$ is a final segment of $I_{l-1}$.  Applying Lemma
\ref{subadd},
$$
D_J((\sigma \otimes \tau)(S)) \leq \sum_{k=l}^L D_J((\sigma
\otimes \tau)(I_k)) + D_J((\sigma \otimes \tau)(S_0)) +
D_J((\sigma \otimes \tau)(S^\prime_0))
$$
By the arguments presented in the proof of the previous
proposition,
\begin{eqnarray*}
D_J((\sigma \otimes \tau)(I_k)) & \leq & \frac{m-1}{m} \\
D_J((\sigma \otimes \tau)(S_0)) & \leq & d(\sigma) + \frac{m-1}{m} \\
D_J((\sigma \otimes \tau)(S^\prime_0)) & \leq & d^\prime(\sigma) +
\frac{m-1}{m}
\end{eqnarray*}
Thus,
$$
D(\sigma \otimes \tau) \leq m - 1 + d(\sigma) + d^\prime(\sigma)
$$
\end{proof}

If we apply these results to a product of permutations,
\begin{cor} \label{prodperm} If, for $i \leq k$, $\sigma_i \in S_{n_i}$, where $n_i > 1$, then
$$
D(\bigotimes_{i=1}^{k} \sigma_i ) \leq n_k + 2 \sum_{i=1}^{k-1}
n_i - 2k + 1.
$$
\end{cor}
\begin{proof}
By Proposition \ref{inh}, $d(\bigotimes_{i=1}^{m} \sigma_i ) \leq
n_m - 1 + d(\bigotimes_{i=1}^{m-1} \sigma_i )$.  Inductively,
then, we find
$$
d(\bigotimes_{i=1}^{m} \sigma_i ) \leq \sum_{i=2}^{m} (n_i - 1) +
d(\sigma_1).
$$
Since $d(\sigma_1) \leq n_1 - 1$, we may write
$d(\bigotimes_{i=1}^{m} \sigma_i ) \leq \sum_{i=1}^{m} n_i - m$
and, similarly, $d^\prime(\bigotimes_{i=1}^{m} \sigma_i ) \leq
\sum_{i=1}^{m} n_i - m$.  Thus,
\begin{align*}
D(\bigotimes_{i=1}^{k} \sigma_i ) & \leq n_k - 1 + d(\bigotimes_{i=1}^{k-1} \sigma_i ) + d^\prime(\bigotimes_{i=1}^{k-1} \sigma_i ) \\
& \leq n_k + 2 \sum_{i=1}^{k-1} n_i - 2k + 1.
\end{align*}
\end{proof}

Corollary \ref{prodperm} provides us with a large family of very
strongly quasirandom permutations.  To see this, let
$\{\sigma_i\}_{i=1}^\infty$ be a sequence of permutations with
$\sigma_i \in S_{n_i}$.  Then, letting $\lambda_k = \sum_{i=1}^{k}
n_i / \prod_{i=1}^{k} n_i$ and applying the corollary,
$$
\frac{D(\bigotimes_{i=1}^{k} \sigma_i)}{\prod_{i=1}^{k} n_i} < 2
\frac{\sum_{i=1}^{k} n_i}{\prod_{i=1}^{k} n_i} = 2 \lambda_k
$$
But,
\begin{align*}
\lambda_k &= \frac{1}{n_k} \cdot \frac{\sum_{i=1}^{k-1} n_i}{\prod_{i=1}^{k-1} n_i} + \frac{n_k}{\prod_{i=1}^{k} n_i} \\
&= \left(\frac{1}{n_k} + \frac{1}{\sum_{i=1}^{k-1} n_i} \right )
\lambda_{k-1}
\end{align*}
Thus, the ratio of the discrepancy to the size of the product
permutations tends to zero quickly.  In particular, $\sigma^{(k)}$
is very strongly quasirandom, since $\sigma^{(k)} \in S_{n^k}$.
That is, if $\sigma^{(k)} \in S_N$, then $D(\sigma^{(k)}) < 2kn =
O(\log N)$.  Immediately one wonders whether permutations exist
with discrepancies which grow slower than $\log N$.  A theorem of Schmidt [\ref{S1}] answers this question in the negative, implying that the $D(\sigma^{(k)})$ are, in a sense, ``maximally'' quasirandom.

\begin{theorem} [Schmidt] \label{schmidt} Let $\{x_i\}_{i=0}^{N-1} \subset [0,1)$, and define
$$
D(m) = \sup_{\alpha \in [0,1)} \left | |\{x_i\}_{i=0}^{m-1} \cap
[0,\alpha) | - m \alpha \right |.
$$
Then there exists an integer $n \leq N$ so that $D(n) > \log N /
100$.
\end{theorem}

We may immediately conclude that discrepancies grow at least as
fast as $\log N$.

\begin{cor} \label{disclowerbound} For any $\sigma \in S_N$, $D(\sigma) > \log N / 100 - 1$.
\end{cor}
\begin{proof}
Take $x_i = \sigma(i)/N$ in Theorem \ref{schmidt}.  Then there
exists an $\alpha \in [0,1)$ and an $n \leq N$ so that
$$
\left | \left |\frac{\sigma([0,n-1])}{N} \cap [0,\alpha) \right |
- n \alpha \right | > \frac{\log N}{100}.
$$
Defining $k=\lfloor \alpha N \rfloor$, we have
$$
\left | |\sigma([0,n-1]) \cap [0,k] | - \frac{n(k+1)}{N} \right |
+ \left | \frac{n(k+1)}{N} - n \alpha \right | > \frac{\log
N}{100}.
$$
Therefore, if we let $I$ and $J$ vary over all intervals in
$\mathbb{Z}_N$,
$$
\max_{I,J} \left | |\sigma(I) \cap J | - \frac{|I||J|}{n} \right |
> \frac{\log N}{100} - 1
$$
so that $D(\sigma) > \log N / 100 - 1$.
\end{proof}

One might expect that the algebraic properties of quasirandom
permutations, such as the number of cycles, should be
approximately that of random permutations (in this case, $\log
n$).  However, we have the following counterexample.  Let $i_n$ be
the identity permutation on $\mathbb{Z}_n$.  Then $i^{(k)}_n$ is
always an involution in $S_{n^k}$ -- and the sequence $\{i^{(k)}_n/n^k\}_{i=0}^{n-1} \subset [0,1)$ is an initial segment of the van der Corput sequence.  In fact, under this interpretation, Corollary \ref{prodperm} can be considered a generalization of the classical theorem that the discrepancy of the van der Corput sequence is $O(\log N)$.  (See, for example, [\ref{F1}] for a modern version of this result.)

\begin{prop} $i^{(k)}_n(x)$ is the element of $\mathbb{Z}_{n^k}$ whose base $n$ expansion is the reverse of the base $n$ expansion of $x$.
\end{prop}
\begin{proof}
The proof is by induction on $k$.  The case $k=1$ is obvious.
Suppose it were true of $i^{(k)}_n$.  Let $x \in
\mathbb{Z}_{n^{k+1}}$ have the base $n$ expansion $x_{k} x_{k-1}
\cdots x_1 x_0$.  Then
\begin{align*}
i^{(k+1)}_n(x) & = (i^{(k)}_n \otimes i_n)(x) \\
& = i_n(\lfloor \frac{x}{n^k} \rfloor) + n i^{(k)}_n(x \! \! \! \! \mod n^k) \\
& = (\lfloor \frac{x}{n^k} \rfloor) + n i^{(k)}_n(x_{k-1} x_{k-2} \cdots x_1 x_0) \\
& = x_k + n (x_0 x_1 \cdots x_{k-2} x_{k-1}) \\
& = x_0 x_1 \cdots x_{k-2} x_{k-1} x_k.
\end{align*}
\end{proof}

\section{{\sc Conclusion}}
\setcounter{theorem}{0}

The original motivation for this paper was a (still unanswered)
question of R. L. Graham [\ref{G2}].  For a sequence of
permutations $\sigma_j \in S_{n_j}$, let {\bf P(k)} be the
property of \textit{asymptotic $k$-symmetry}: for each $\tau \in
S_k$,
$$
\left | X^{\tau}(\sigma_j) - \frac{\binom{n_j}{k}}{k!} \right | =
o(n_j^k).
$$
Note that this property is weaker than property \textbf{[kS]} of
Theorem \ref{qrp}, which we will call \textit{strong} asymptotic
$k$-symmetry.  Theorem \ref{qrp} says that strong asymptotic
$k$-symmetry implies strong asymptotic $(k+1)$-symmetry for any $k
\geq 2$.  Graham asks whether there exists an analogous $N$ so
that, for all $k > N$, $\mbox{\bf P(k)} \Rightarrow \mbox{\bf
P(k+1)}$?  At first it might seem like one is asking for too much.
However, precisely this type of phenomenon occurs for graphs
([\ref{CGW1}]).  It turns out that, if we let {\bf G(k)} be the
property that all graphs on $k$ vertices occur as subgraphs at
approximately the same rate, then
$$
\mbox{\bf G(1)} \Leftarrow \mbox{\bf G(2)} \Leftarrow \mbox{\bf
G(3)} \Leftarrow \mbox{\bf G(4)} \Leftrightarrow \mbox{\bf G(5)}
\Leftrightarrow \mbox{\bf G(6)} \Leftrightarrow \cdots
$$
In particular, {\bf G(4)} implies quasirandomness, which in turn
implies {\bf G(k)} for all $k$.

The fact that $\mbox{\bf P(1)} \not \Rightarrow \mbox{\bf P(2)}$
is trivial.  To show that $\mbox{\bf P(2)} \not \Rightarrow
\mbox{\bf P(3)}$, let $\sigma_n \in S_{2n}$ be the permutation
which sends $x$ to $x+n$.  Then $X^{01}(\sigma_n) = 2n(2n-1)$, and
$X^{10}(\sigma_n) = 4n^2$, so that $\left | X^{01}(\sigma_n) -
X^{10}(\sigma_n) \right | = o((2n)^2)$.  However, the pattern
$(021)$ {\it never} appears in $\sigma_n$.  We have been unable to
date to provide an analogous result for any $P(k)$ with $k
> 2$.

A second, very natural question is that of the existence of
\textit{perfect} $m$-symmetry: the property of having all
subsequence statistics \textit{precisely} equal to their expected
values.  That is, for $\sigma \in S_n$,
$$
X^{\tau}(\sigma) = \frac{\binom{n}{k}}{k!}
$$
for all $\tau \in S_m$.  For this to occur, the number of
permutations of length $m$ must evenly divide $\binom{n}{m}$.  Let
\textbf{D($m$)} be the property of an integer $N$ that
$$
m! \left | \binom{n}{m} \right. \!\! .
$$
It is easy to see that a permutation $\sigma \in S_n$ with perfect
$m$-symmetry must have perfect $m^\prime$-symmetry for any
$m^\prime \leq m$, so $n$ must satisfy \textbf{D($m^\prime$)} for
all such $m^\prime$.  Let $h(m)$ be the least $n$ for which this
occurs.  A quick calculation reveals that $h(2)=4$, $h(3)=9$,
$h(4)=64$, and $h(5)=128$.  In fact, there is a perfect
$2$-symmetric permutation on $4$ symbols: 3012.  A computer search
revealed that there are exactly two $3$-symmetric permutations on
$9$ symbols: 650147832 and its reverse, 238741056.  No
$m$-symmetric permutation is known for $m>3$, and the question of
whether such permutations exist remains open.  We conjecture that
an $m$-symmetric permutation on sufficiently many symbols exists
for all $m$, and believe it likely that one exists on $h(m)$
symbols.

Finally, the selection of intervals as the sets which measure ``proximity''
in the definition of quasirandomness was a natural but somewhat arbitrary choice.
It would be worth investigating the properties of ``$(\mathcal{A},\mathcal{B})$-quasirandom''
permutations for families $\mathcal{A}$, $\mathcal{B} \subset 2^{\mathbb{Z}_n}$, i.e., permutations $\sigma$ such that $\max_{A,B} D_B(\sigma(A)) = o(n)$ for $A \in \mathcal{A}$ and $B \in \mathcal{B}$.

\section*{{\sc Acknowledgements}}

The author wishes to thank Fan Chung Graham and Ron Graham for
their tremendous help in formulating and attacking the problems discussed above.  He also thanks Chris Dillard, Robert Ellis, and Lei Wu for helpful discussions during the development of this work.

\section*{{\sc References}}
\begin{enumerate}
\item \label{Ah1} {\sc L. V. Ahlfors}, ``Complex Analysis,'' McGraw-Hill Book Co., New York, 1978.
\item \label{A1} {\sc N. Alon}, Generating Pseudo-Random Permutations and Maximum Flow Algorithms, {\it Inform. Process. Lett.} {\bf 35} (1990), 201--204.
\item \label{AF1} {\sc N. Alon and E. Friedgut}, On the number of permutations avoiding a given pattern, {\it J. Comb. Theory Ser. A} {\bf 89} (2000), 133--140.
\item \label{AS1} {\sc N. Alon and J. H. Spencer}, ``The probabilistic method,'' Wiley-Interscience Series in Discrete Mathematics and Optimization. Wiley-Interscience [John Wiley \& Sons], New York, 2000. 
\item \label{B1} {\sc M. B\'{o}na}, The solution of a conjecture of Stanley and Wilf for all layered patterns, {\it J. Comb. Theory Ser. A} {\bf 85} (1999), 96--104.
\item \label{Ch1} {\sc B. Chazelle}, ``The Discrepancy Method,'' Cambridge University Press, Cambridge, 2000.
\item \label{CG1} {\sc F. R. K. Chung and R. L. Graham}, Quasi-random set systems, {\it J. Amer. Math. Soc.} {\bf 4} (1991), 151--196.
\item \label{CG2} {\sc F. R. K. Chung and R. L. Graham}, Quasi-random subsets of $Z\sb n$, {\it J. Combin. Theory Ser. A} {\bf 61} (1992), 64--86.
\item \label{CGW1} {\sc F. R. K. Chung, R. L. Graham, and R. M. Wilson}, Quasi-random graphs, {\it Combinatorica} {\bf 9} (1989), 345--362.
\item \label{D1} {\sc P. J. Davis}, ``Circulant Matrices,'' Wiley, New York, 1979.
\item \label{D2} {\sc J. A. Davis, R. C. Entringer, R. L. Graham, and G. J. Simmons}, On permutations containing no long arithmetic progressions, {\it Acta Arithmetica}, {\bf 34} (1977/78), 81--90.
\item \label{ES1} {\sc P. Erd\H{o}s and G. Szekeres}, A combinatorial problem in geometry, {\it Compocito Math.} {\bf 2} (1935), 464--470.
\item \label{F1} {\sc H. Faure}, Discr\'{e}pance quadratique de suites infinies en dimension un. Th\'{e}orie des nombres (Quebec, PQ, 1987), 207--212, de Gruyter, Berlin, 1989. 
\item \label{G1} {\sc W. T. Gowers}, A new proof of Szemer\'edi's Theorem, {\it Geometric and Functional Analysis} {\bf 11} (2001), 465--588.
\item \label{G2} {\sc R. L. Graham}, personal communication.
\item \label{L1} {\sc J. C. Lagarias}, Well-spaced labelling of points in rectangular grids, {\it SIAM J. Discrete Math.} {\bf 13} (2000), 521--534.
\item \label{N1} {\sc H. Niederreiter}, Quasi-Monte Carlo methods and pseudo-random numbers. {\it Bull. Amer. Math. Soc.} {\bf 84} (1978), 957--1041. 
\item \label{MMW1} {\sc B. D. McKay, J. Morse, and H. S. Wilf}, The distributions of the entries of Young tableaux, to appear.
\item \label{Sa1} {\sc V. N. Sachkov}, ``Probabilistic Methods in Combinatorial Analysis,'' Encyclopedia of Mathematics and its Applications {\bf 56}, Cambridge University Press, Cambridge, 1997.
\item \label{SS1} {\sc F. W. Schmidt and R. Simion}, Restricted permutations, {\it European J. Combin.} {\bf 6} (1985), 383--406.
\item \label{S1} {\sc W. M. Schmidt}, Irregularities of distribution VII, {\it Acta Arith.} {\bf 21} (1972), 45–-50.
\item \label{St1} {\sc R. P. Stanley}, Log-concave and unimodal sequences in algebra, combinatorics, and geometry. Graph theory and its applications: East and West (Jinan, 1986), 500--535, Ann. New York Acad. Sci., 576, New York Acad. Sci., New York, 1989.
\item \label{St2} {\sc R. P. Stanley}, Enumerative combinatorics, Vol. 1, {\it Cambridge Studies in Advanced Mathematics} {\bf 49}, Cambridge University Press, Cambridge, 1997.
\end{enumerate}

\ifx\ispreprint\undefined
\pagebreak
\def\baselinestretch{1}
\begin{center}
\begin{figure}[ht]
\hspace{.7in} \input{figa.pstex_t} \caption{Diagram of
implications for Theorem \ref{balance}.}
\end{figure}
\end{center}
\else
\fi

\ifx\ispreprint\undefined
\pagebreak
\def\baselinestretch{1}
\begin{center}
\begin{figure}[ht]
\hspace{1.25in} \input{figb.pstex_t} \caption{Diagram of
implications for Theorem \ref{qrp}.}
\end{figure}
\end{center}
\else
\fi

\end{document}